\documentclass[12pt,letterpaper]{article}

\usepackage{amssymb}
\usepackage{amsmath}
\usepackage{amsthm}
\usepackage{graphicx}
\usepackage{graphics}
\usepackage{setspace}
\usepackage{subfigure}
\usepackage[margin=1in]{geometry}

\begin{document}

\newcommand{\pardfdx}{$\frac{\partial \mathbf{f}}{\partial \mathbf{x}} $}
\newcommand{\pardfdxnm}{\frac{\partial \mathbf{f}}{\partial \mathbf{x}}}
\newcommand{\pardfdynm}{\frac{\partial \mathbf{f}}{\partial \mathbf{y}}}
\newcommand{\boldf}{\mathbf{f}}
\newcommand{\boldg}{\mathbf{g}}
\newcommand{\boldx}{\mathbf{x}}
\newcommand{\boldz}{\mathbf{z}}
\newcommand{\boldT}{\mathbf{T}}
\newcommand{\boldM}{\mathbf{M}}
\newcommand{\boldF}{\mathbf{F}}
\newcommand{\boldR}{\mathbf{R}}
\newcommand{\boldy}{\mathbf{y}}
\newcommand{\boldu}{\mathbf{u}}
\newcommand{\pardfdxdel}{\frac{\partial \boldf_{\delta}}{\partial \boldx}}
\newcommand{\pardfdxgam}{\frac{\partial \boldf_{\gamma}}{\partial \boldx}}

\newcommand{\R}{{\mathbb R}}

\newtheorem{definition}{Definition}
\newtheorem{theorem}{Theorem}
\newtheorem{lemma}{Lemma}

\newcommand{\todo}[1]{\vspace{5 mm}\par \noindent \marginpar{\textsc{ToDo}}
\framebox{\begin{minipage}[c]{0.95 \columnwidth} \tt #1
\end{minipage}}\vspace{5 mm}\par}

\title{Stability and Robustness Analysis of Nonlinear Systems \\
via Contraction Metrics and SOS Programming }

\author{Erin M. Aylward$^1$ \qquad Pablo A. Parrilo$^1$ \qquad Jean-Jacques E. Slotine$^2$ \\
$^1$Laboratory for Information and Decision Systems\\
$^2$Nonlinear Systems Laboratory \\
 Massachusetts Institute of Technology \\
Cambridge, MA 02139, USA}

\date{LIDS Technical Report \#2691 \\ \vspace{.6cm} \today}

\maketitle

 \tableofcontents
 \pagebreak
\begin{abstract}

Contraction analysis is a stability theory for nonlinear systems where
stability is defined incrementally between two arbitrary
trajectories. It provides an alternative framework in which to study
uncertain interconnections or systems with external inputs, where it
offers several significant advantages when compared with traditional
Lyapunov analysis. Contraction-based methods are particularly useful
for analyzing systems with uncertain parameters and for proving
synchronization properties of nonlinear oscillators. Existence of a
contraction metric for a given system is a necessary and sufficient
condition for global exponential convergence of system
trajectories. For systems with polynomial or rational dynamics, the
search for contraction metrics can be made fully algorithmic through
the use of convex optimization and sum of squares (SOS)
programming. The search process is made computationally tractable by
relaxing matrix definiteness constraints, whose feasibility indicate
existence of a contraction metric, into SOS constraints on polynomial
matrices.  We illustrate the results through examples from the
literature, emphasizing the advantages and contrasting the differences
between the contraction approach and traditional Lyapunov techniques.

\end{abstract}

\section{Introduction}
Contraction analysis is a stability theory for nonlinear systems where
stability is defined incrementally between two arbitrary trajectories
\cite{OCA}. The existence of a contraction metric for a nonlinear
system ensures that a suitably defined distance between nearby
trajectories is always decreasing, and thus trajectories converge
exponentially and globally. One important application of contraction
theory is its use in studying the synchronization of nonlinear coupled
oscillators \cite{OPC-NO}.  These oscillators present themselves in a
variety of research fields such as mathematics, biology, neuroscience,
electronics, and robotics. The use of coupled oscillators in each of
these fields, as well as how contraction theory can be used to analyze
networks of coupled identical nonlinear oscillators can be can be
found in \cite{OPC-NO} and the references listed therein.

Contraction theory nicely complements Lyapunov theory, a standard
nonlinear stability analysis technique, as it provides an
alternative framework in which to study convergence and robustness
properties of nonlinear systems. For autonomous systems one can
interpret the search for a contraction metric as the search for a
Lyapunov function with a certain structure. This statement will be
explained further in Section \ref{UD}.  There are, however,
advantages to searching for a contraction metric instead of
searching explicitly for a Lyapunov function. In particular, as we
will show, contraction metrics are useful for analyzing uncertain
nonlinear systems. In general, nonlinear systems with uncertain
parameters can prove quite troublesome for standard Lyapunov
methods, since the uncertainty can change the equilibrium point of
the system in very complicated ways, thus forcing the use of
parameter-dependent Lyapunov functions in order to prove stability
for a range of the uncertain parameter values.

Much of the literature on parameter-dependent Lyapunov functions
focuses on linear systems with parametric uncertainty
\cite{PDL:Haddad,PDL:Gahinet,PDL:Feron, Rantzerpaper}. However, if
a linear model is being used to study a nonlinear system around an
equilibrium point, changing the equilibrium of the nonlinear
system, necessitates relinearization around the new equilibrium.
If the actual position of the equilibrium, in addition to the
stability properties of the equilibrium, of the nonlinear system
depends on the uncertainty, it may be impossible to obtain any
kind of closed form expression of the equilibrium in terms of the
uncertain parameters. Thus, parameterizing the linearization in
terms of the uncertainty may not be an option.

A well-studied method of dealing with specific forms of nonlinearities
is to model the nonlinear system as a linear system with bounded
uncertainty. In particular, in \cite{LMIBook} polytopic linear
differential inclusions (LDIs), norm-bound LDIs, and diagonal
norm-bound LDIs are considered. These techniques are computationally
tractable as they reduce to convex optimization problems. Though these
methods work for various kinds of uncertainty, it is also desirable to
find methods to study the stability of nonlinear systems that do not
easily admit linear approximations with the nonlinearities covered
with uncertainty bounds.

Contraction theory provides a framework in which to study the
stability behavior of more general uncertain nonlinear systems.
This framework eliminates many of the restrictions and problems
that may be encountered when trying to analyze uncertain nonlinear
systems with traditional linearization techniques or Lyapunov
methods. This results from the fact that if a nominal system is
contracting with respect to a certain contraction metric, it is
often the case that the uncertain system with additive or
multiplicative uncertainty within a certain range will still be
contracting with respect to the same metric, even if the
perturbation changes the position of the equilibrium of the
system. Thus, it is possible to determine stability of the system
for a range of values of the uncertain parameter without
explicitly tracking how the uncertainty changes the location of
the equilibrium. These ideas will be discussed further in
Section~\ref{UD}.

Another interesting feature of the contraction framework is its
relative flexibility in incorporating inputs and outputs. For
instance, to prove contraction of a class of systems with external
inputs, it is sufficient to show the existence of a contraction
metric with a certain structure.  This feature, which will be
discussed in Section~\ref{EI}, is central in using contraction
theory to prove synchronization of coupled nonlinear oscillators.

To translate the theoretical discussion above into effective
practical tools, it is desirable to have efficient computational
methods to numerically obtain contraction metrics. Sum of squares
(SOS) programming provides one such method. SOS programming is
based on techniques that combine elements of computational algebra
and convex optimization, and has been recently used to provide
efficient convex relaxations for several computationally hard
problems \cite{CAND}. In this paper we will show how SOS
programming enables the search for contraction metrics for the
class of nonlinear systems with polynomial dynamics. We discuss
how to use SOS methods to find bounds on the maximum amount of
uncertainty allowed in a system in order for the system to retain
the property of being contracting with respect to the contraction
metric of the unperturbed system. We also use SOS methods to
optimize the contraction matrix search to obtain a metric that
provides the largest symmetric uncertainty interval for which we
can prove the system is contracting.

This paper is organized as follows: in Section~\ref{CT} we give
background material on contraction theory. Section~\ref{SOSSect}
discusses sum of squares (SOS) polynomials and matrices.  We
present next an algorithm which uses SOS programming to
computationally search for contraction metrics for nonlinear
systems. We discuss why contraction theory is useful for studying
systems with uncertain dynamics in Section~\ref{UD} and external
inputs in Section~\ref{EI}. Finally, in Section~\ref{Conclusions}
we present our conclusions, and outline possible directions for
future work.

\section{Contraction Analysis}\label{CT}

Contraction analysis is a relatively recently developed stability
theory for nonlinear systems analysis \cite{OCA}. The theory attempts
to answer the question of whether the limiting behavior of a given
dynamical system is independent of its initial conditions.  More
specifically, contraction analysis is a theory in which stability is
defined incrementally between two arbitrary trajectories.  It is used
to determine whether nearby trajectories converge to one another. This
section summarizes the main elements of contraction analysis; a much
more detailed account can be found in \cite{OCA}.

We consider deterministic dynamical systems of the form
\begin{equation}
\dot{\boldx} = \boldf(\boldx(t),t), \label{sysdynamics}
\end{equation}
where $\boldf$ is a nonlinear vector field and $\boldx (t)$ is an
$n$-dimensional state vector. For this analysis it is assumed that
all quantities are real and smooth and thus that all required
derivatives or partial derivatives exist and are continuous. This
existence and continuity assumption clearly holds for polynomial
vector fields.

Under the assumption that all quantities are real and smooth, from
equation (\ref{sysdynamics}) we can obtain the differential
relation
\begin{equation}
\delta \dot{\boldx} (t) = \frac{\partial \boldf}{\partial \boldx
(t)}(\boldx (t),t)\delta \boldx (t), \label{diffdynamics}
\end{equation}
where $\delta \boldx (t)$ is an infinitesimal displacement at a
fixed time. For notational convenience from here on we will write
$\boldx$ for $\boldx(t)$, but in all calculations it should be
noted that $\boldx$ is a function of time.

The infinitesimal squared distance between two trajectories is
$\delta \boldx^T \delta \boldx$.  Using (\ref{diffdynamics}), the
following equation for the rate of change of the squared distance
between two trajectories is obtained:
\begin{equation} \frac{d}{dt}(\delta \boldx^T \delta \boldx) = 2\delta \boldx^T
\delta \dot{\boldx} =  2 \delta \boldx^T \frac{\partial
\boldf}{\partial \boldx} \delta \boldx. \label{eq:diffdynamics2}
\end{equation} If $\lambda_1(\boldx,t)$ is the largest eigenvalue of the
symmetric part of the Jacobian $\frac{\partial \boldf}{\partial
\boldx}$ (i.e. the largest eigenvalue of $\frac{1}{2}(\pardfdxnm +
\pardfdxnm^T)$),  then it follows from (\ref{eq:diffdynamics2}) that
\begin{equation}
\frac{d}{dt}(\delta \boldx^T \delta \boldx) \leq 2
\lambda_1(\boldx, t) \delta \boldx ^T \delta \boldx.
\end{equation} Integrating both sides gives
\begin{equation}
||\delta \boldx|| \leq || \delta \boldx _o|| \, e^{\int_0^t
\lambda_1(\boldx,t)dt}.
\label{expdyn}
\end{equation} If $\lambda_1(\boldx,t)$ is uniformly strictly negative
(i.e. $ ( \frac{\partial \boldf}{\partial \boldx} + \frac{\partial
\boldf}{\partial \boldx}^T ) \prec 0 \; \forall \; \boldx, t)$, it
follows from (\ref{expdyn}) that any infinitesimal length
$||\delta \boldx||$ converges exponentially to zero.  By path
integration the distance of any finite path also converges
exponentially to zero.

A more general definition of length can be given by
\begin{equation}
\delta \boldz^T \delta \boldz =
\delta \boldx ^T \boldM (\bold x, t) \delta \boldx
\label{RiemannSp}
\end{equation}
where $\boldM(\boldx, t)$ is a symmetric, uniformly positive
definite and continuously differentiable metric (formally, this
defines a Riemannian manifold).  This notion of infinitesimal
distance defined with respect to a metric can be use to define a
finite distance measure between two trajectories with respect to
this metric. Specifically, the distance between two points $P_1$
and $P_2$ with respect to the metric $\mathbf{M}(\boldx,t)$ is
defined as the shortest path length, in other words the smallest
path integral $\int_{P_1}^{P_2} \sqrt{\delta \boldx^T
\boldM(\boldx,t) \delta \boldx}$, between these two points.
Accordingly a ball of center $\mathbf{c}$ with radius $R$ is
defined as the set of all points whose distance to $\mathbf{c}$
with respect to $\boldM (\boldx,t)$ is strictly less than $R$.

Under the definition of infinitesimal length given in
(\ref{RiemannSp}), the equation for its rate of change becomes
\begin{equation}\label{contractcond2} \frac{d}{dt}(\delta \boldx
^T \boldM \delta \boldx) = \delta \boldx^T(\pardfdxnm ^T \boldM +
\boldM
\pardfdxnm + \dot{\boldM}) \delta \boldx \end{equation}
where $\boldM$ is shorthand notation for $\boldM(\boldx,t)$.
Convergence to a single trajectory occurs in regions where
$(\pardfdxnm ^T \boldM + \boldM
\pardfdxnm + \dot{\boldM})$ is
uniformly negative definite.   It should be noted that
$\dot{\boldM}  = \dot{\boldM}(\boldx, t)= \frac{\partial
\boldM(\boldx,t)}{\partial \boldx}\frac{d \boldx}{dt} +
\frac{\partial \boldM (\boldx, t)}{\partial t}$. The above
analysis leads to the following definition and theorem:

\begin{definition}[\cite{OCA}]
Given the system equations $\dot{\boldx} = \boldf(\boldx,t)$, a
region of the state space is called a \emph{contraction region}
with respect to a uniformly positive definite metric
$\boldM(\boldx, t)$  if $(\pardfdxnm ^T \boldM + \boldM
\pardfdxnm + \dot{\boldM})$ is uniformly negative definite in
that region.
\end{definition}

\begin{theorem}[\cite{OCA}]
Consider the system equations $\dot{\boldx} = \boldf(\boldx,t)$.
Assume a trajectory starts in a ball of constant radius that is
defined with respect to the metric $\boldM(\boldx, t)$, that is
centered at a another given trajectory, and that is contained at
all times in a contraction region with respect to the metric
$\boldM(\boldx,t)$. Then the first trajectory will remain in that
ball and converge exponentially to the center trajectory.
Furthermore, global exponential convergence to the center
trajectory is guaranteed if the whole state space is a contraction
region with respect to the metric $\boldM(\boldx,t)$.
\end{theorem}

Definition 1 provides sufficient conditions for a system to be
contracting. Namely, the following should be satisfied:
\begin{enumerate}
\item The matrix $\boldM(\boldx, t)$ must be a uniformly positive
definite matrix, i.e.,
\begin{equation}\label{Mpos}\boldM(\boldx, t)
\succeq \epsilon \mathbf{I} \succ 0 \qquad \forall \boldx, t.
\end{equation}

\item The metric variation $\pardfdxnm ^T \boldM + \boldM
\pardfdxnm + \dot{\boldM}$ must be a uniformly negative definite
matrix, i.e.,
\begin{equation}\label{Rneg}\boldR(\boldx,t) = \pardfdxnm ^T \boldM + \boldM
\pardfdxnm + \dot{\boldM} \preceq -\epsilon \mathbf{I} \prec 0 \; \; \; \forall \boldx,t.
\end{equation}
\end{enumerate}An explicit rate of convergence of trajectories $\beta$ can be
found by finding a $\boldM(\boldx, t)$ that satisfies (\ref{Mpos})
and
\begin{equation}
\pardfdxnm ^T \boldM + \boldM
\pardfdxnm + \dot{\boldM} \preceq -\beta \boldM .
\end{equation}

The notation above is standard; $\succ$, and $\succeq$ mean
positive definite and positive semidefinite respectively, while
$\prec$ and $\preceq$ mean negative definite and negative
semidefinite respectively. If the system dynamics are linear and
$\boldM(\boldx,t)$ is constant (i.e. $\boldM(\boldx,t) = M$), the
conditions above reduce to those in standard Lyapunov analysis
techniques. Lyapunov theory shows that the system $\dot{\boldx}(t)
= A\boldx(t)$ is stable (i.e., all trajectories converge to 0) if
and only if there exists a positive definite matrix $M$ (i.e., $M
\succ 0$) such that $A^TM + MA \prec 0$.

It should be noted that if a global contraction metric exists for
an autonomous system, all trajectories converge to a unique
equilibrium point, and we can always produce a Lyapunov function
for the system from the contraction metric \cite{OCA}. We assume,
without loss of generality, that the equilibrium is at the origin.
If the system dynamics are $\boldf(\boldx)$ and $\boldM(\boldx)$
is a time-invariant contraction metric for the system, then
$V(\boldx) = \boldf(\boldx)^T \boldM(\boldx)\boldf(\boldx)$ is a
Lyapunov function for the system since $V(\boldx) > 0$ and
$\dot{V} = \boldf(\boldx)^T(\frac{\partial \boldf}{\partial
\boldx}^T \boldM + \bold M \frac{\partial \boldf}{\partial \boldx}
+ \dot{\boldM})\boldf(\boldx)\leq -\beta V$. This shows that
$\dot{\boldx} = \boldf(\boldx)$ tends to $\mathbf{0}$
exponentially, and thus that $\boldx$ tends towards a finite
equilibrium point.

For a constant metric $\boldM(\boldx,t) = M$, this reduces to
Krasovskii's Method \cite{Khalil}. We note that for systems with
uncertainty there are good reasons to search for a contraction
metric to create Lyapunov function of this structure instead of
searching for a Lyapunov function directly. These reasons will
become clear in Section \ref{UD}.

The problem of searching for a contraction metric thus reduces to
finding a matrix function $\boldM(\boldx,t)$ that satisfies the
conditions above. As we will see, SOS methods will provide a
computationally convenient approach to this problem.

\section{Sum of Squares (SOS) Polynomials and Programs}
\label{SOSSect}

The main computational difficulty of problems involving constraints
such as the ones in (\ref{Mpos}) and (\ref{Rneg}) is the lack of
efficient numerical methods that can effectively handle multivariate
nonnegativity conditions. A convenient approach for this, originally
introduced in \cite{Parrilo:Phd}, is the use of sum of squares (SOS)
relaxations as a suitable replacement for nonnegativity. We present
below the basic elements of these techniques.

A multivariate polynomial $p(x_1, x_2, ..., x_n) = p(\boldx) \in
\R[\boldx]$ is a sum of squares (SOS) if there exist polynomials
$f_1(\boldx), ..., f_m(\boldx) \in \R[\boldx]$ such that
\begin{equation}\label{SOSeqn}
p(\boldx) = \sum_{i=1}^m f_i^2(\boldx).
\end{equation}
The existence of a SOS representation for a given polynomial is a
sufficient condition for its global nonnegativity, i.e., equation
(\ref{SOSeqn}) implies that $p(\boldx) \geq 0 \; \forall \; \boldx
\in \R^n$. The SOS condition (\ref{SOSeqn}) can be shown to be
equivalent to the existence of a positive semidefinite matrix $Q$
such that
\begin{equation}\label{Q}
p(\boldx) = Z^T(\boldx)QZ(\boldx)
\end{equation} where $Z(\boldx)$ is a vector of monomials of
degree less than or equal to deg($p$)/2.  This equivalence of
descriptions between (\ref{SOSeqn}) and (\ref{Q}) makes finding an
SOS decomposition a computationally tractable procedure. Finding a
symmetric positive semidefinite $Q$ subject to the affine
constraint (\ref{Q}) is a semidefinite programming problem
\cite{Parrilo:Phd,SOSTOOLS}.

Using the notion of a SOS polynomial as a primitive, we can now
introduce a convenient class of optimization problems.  A
\emph{sum of squares program} is a convex optimization problem of
the form:
\begin{eqnarray}
& \min& \sum_{j = 1}^J w_j \, c_j \nonumber \\
& \text{subject to}&\; \; a_{i,0} + \sum_{j =1}^J a_{i,j}(\boldx)
\, c_j \quad \mbox{ is SOS for }\, i= 1, ..., I, \nonumber
\end{eqnarray}
where the $c_j$'s are the scalar real decision variables, the $w_j$'s
are given real numbers that define the objective function, and the
$a_{i,j}(\boldx)$ are given multivariate polynomials.  There has
recently been much interest in SOS programming and SOS optimization as
these techniques provide convex relaxations for various
computationally hard optimization and control problems; see e.g.
\cite{Parrilo:Phd,sdprelax,Lasserre,CAND} and the volume \cite{GarulliHenrion}.

A SOS decomposition provides an explicit certificate of the
nonnegativity of a scalar polynomial for all values of the
indeterminates. In order to design an algorithmic procedure to search
for contraction metrics, we need to introduce a similar idea to ensure
that a polynomial matrix is positive definite for every value of the
indeterminates. A natural definition is as follows:

\begin{definition}[\cite{ParriloGatermann}]
Consider a symmetric matrix with polynomial entries
$\mathbf{S}(\boldx) \in \R[\boldx]^{m \times m}$, and let $\boldy =
[y_1, \ldots, y_m]^T$ be a vector of new indeterminates. Then
$\mathbf{S}(\boldx) $ is a \emph{sum of squares matrix} if the scalar
polynomial $\boldy^T \mathbf{S}(\boldx) \boldy$ is a sum of squares in
$\R[\boldx, \boldy]$.
\end{definition} \noindent For notational convenience, we also define a stricter notion:
\begin{definition}
A matrix $\mathbf{S}(\boldx)$ is \emph{strictly SOS} if
$\mathbf{S}(\boldx) - \epsilon \mathbf{I}$ is a SOS matrix for some
$\epsilon > 0$.
\end{definition}

Thus, a strictly SOS matrix is a matrix with polynomial entries
that is positive definite for every value of the indeterminates.
An equivalent definition of an SOS matrix can be given in terms of
the existence of a polynomial factorization: $\mathbf{S}(\boldx)$
is a SOS matrix if and only if it can be decomposed as
$\mathbf{S}(\boldx) = \mathbf{T}(\boldx)^T \mathbf{T}(\boldx)$
where $\mathbf{T}(\boldx) \in \R[\boldx]^{p \times m}$.  For
example,
\begin{displaymath}
\boldM(\boldx) = \left[\begin{array}{cc}
  \omega^2 + \alpha^2(x^2 + k)^2    & \alpha(x^2 + k) \\
  \alpha(x^2 + k ) & 1 \\
\end{array}\right]\end{displaymath} is a SOS matrix for all values of $\alpha$ and $k$.
Indeed, this follows from the decomposition $\boldM(\boldx) =
\mathbf{T}(\boldx)^T \mathbf{T}(\boldx)$, where
\begin{displaymath} \mathbf{T}(\boldx) =
\left[\begin{array}{cc}
  \omega & 0 \\
  \alpha(x^2 + k) & 1 \\
\end{array} \right].\\
\end{displaymath}
SOS matrices have also been used recently by Hol and Scherer
\cite{HolScherer} and Kojima \cite{KojimaSOS} to produce
relaxations of polynomial optimization problems with matrix
positivity definiteness constraints.

\section{Computational Search for Contraction Metrics via SOS Programming}

As explained in Section~\ref{CT}, given a dynamical system, the
conditions for a contraction metric to exist in regions of the
state-space are given by a pair of matrix inequalities. In the
case of metrics $\boldM(\boldx)$ that do not depend explicitly on
time, relaxing the matrix definiteness conditions in (\ref{Mpos})
and (\ref{Rneg}) to SOS matrix based tests makes the search for
contracting metrics a computationally tractable procedure. More
specifically, the matrix definiteness constraints on
$\boldM(\boldx)$ (and $\boldR(\boldx)$) can be relaxed to SOS
matrix constraints by changing the inequality $\boldM(\boldx) -
\epsilon \mathbf{I} \succeq 0$ in~(\ref{Mpos}) (where $\epsilon$
is an arbitrarily small constant) to the weaker condition that
$\boldM(\boldx)$ be a strictly SOS matrix.  With these
manipulations we see that existence of SOS matrices
$\boldM(\boldx)$ and $\boldR(\boldx)$ is a sufficient condition
for contraction.

\begin{lemma}
Existence of a strictly SOS matrix $\boldM(\boldx)$ and a strictly SOS
matrix $-\boldR(\boldx) = -(\pardfdxnm ^T \boldM + \boldM\pardfdxnm +
\dot{\boldM})$ is a sufficient condition for global contraction of an
autonomous system $\dot{\boldx} = \boldf(\boldx)$ with polynomial
dynamics.
\end{lemma}

\begin{proof}
By Theorem 1, a sufficient condition for contraction of any nonlinear
system is the existence of uniformly positive definite
$\boldM(\boldx)$ and $-\boldR(\boldx)$. A sufficient condition for
uniform positive definiteness of $\boldM(\boldx)$ and
$-\boldR(\boldx)$ is the existence of strictly SOS matrices
$\boldM(\boldx)$ and $-\boldR(\boldx)$.
\end{proof}

This lemma can easily be extended to existence of certain SOS
matrices implying contraction with a convergence rate $\beta$ by
redefining $\boldR(\boldx)$ as $\boldR (\boldx) =
\pardfdxnm ^T \boldM + \boldM\pardfdxnm + \dot{\boldM} + \beta
\boldM$. At this point, we do not know if the full converse of
Lemma~1 holds. If a system is exponentially contracting, it is
known that a contraction matrix always exists \cite{OCA}.
Nevertheless, a system with polynomial dynamics may certainly be
contracting under non-polynomial metrics. Furthermore, even if a
positive definite contraction matrix with polynomial entries
$\bold M$ exists, it may not be the case that it is a SOS matrix.
We notice, however, that some of these issues, such as the gap
between ``true'' contracting metrics and SOS-based ones, can be
bridged by using the more advanced techniques explained in
\cite{sdprelax}.

\subsection{Search Algorithm}
\label{SearchAlgorithm}

One main contribution of this work is to show how sum of squares
(SOS) techniques can be used to algorithmically search for a
time-invariant contraction metric for nonlinear systems with
polynomial dynamics. Existence of a contraction metric for
nonlinear systems certifies contraction (or convergence) of system
trajectories. For systems with polynomial dynamics, we can obtain
a computationally tractable search procedure by restricting
ourselves to a large class of SOS-based metrics.

As suggested by Lemma~1, the main idea is to relax the search for
matrices that satisfy matrix definiteness constraints $\boldM(\boldx)
\succ 0$ and $-\boldR(\boldx) \succ 0$ into SOS-matrix sufficient
conditions. Equivalently, we want to find a polynomial matrix
$\boldM(\boldx)$ that satisfies SOS matrix constraints on
$\boldM(\boldx)$ and $\boldR(\boldx)$.
The SOS feasibility problem can then be formulated as finding
$\boldM(\boldx)$ and $\boldR(\boldx)$ such that
$\boldy^T\boldM(\boldx)\boldy$ is SOS and
$-\boldy^T\boldR(\boldx)\boldy$ is SOS.

More specifically, the detailed steps in the algorithmic search of
contraction metrics for systems with polynomial dynamics are as
follows: \\

\noindent 1. Choose the degree of the polynomials in the contraction
metric, and write an affine parametrization of the symmetric matrices
of that degree. For instance, if the degree is equal to two, the
general form of $\boldM(\boldx)$ is

{\small \[ \left[
\begin{array}{cc}
  a_1x_1^2 + a_2x_1x_2 + a_3x_2^2 + a_4x_1 + a_5x_2 + a_6 &  b_1x_1^2 + b_2x_1x_2 + b_3x_2^2 + b_4x_1 + b_5x_2 + b_6 \\
  b_1x_1^2 + b_2x_1x_2 + b_3x_2^2 + b_4x_1 + b_5x_2 + b_6 & c_1x_1^2 + c_2x_1x_2 + c_3x_2^2 + c_4x_1 + c_5x_2 + c_6 \\
\end{array}\right] \]}%
where $a_i$, $b_i$, and $c_i$ are
unknown coefficients.\\

\noindent 2. Calculate $\pardfdxnm$ and define $\boldR(\boldx) :=
\pardfdxnm ^T \boldM + \boldM \pardfdxnm + \dot{\boldM}$. Thus,
$\boldR(\boldx)$ will also be a symmetric matrix with entries that
depend affinely on the same unknown coefficients $a_i$, $b_i$, and
$c_i$. \\

\noindent 3. Change matrix constraints $\boldM(\boldx) \succ 0$
$\forall \boldx$, and $\boldR(\boldx) = \pardfdxnm ^T \boldM +
\boldM
\pardfdxnm + \dot{\boldM} \prec 0 \; \forall \: \boldx$ into scalar
constraints on quadratic functions $p(\boldx,\boldy)
=\boldy^T\boldM(\boldx)\boldy > 0 \; \forall \: \boldx , \:
\boldy$, and $r(\boldx,\boldy) = \boldy^T\boldR(\boldx)\boldy =
\boldy^T(\pardfdxnm ^T \boldM + \boldM
\pardfdxnm + \dot{\boldM})\boldy <0 \; \forall \: \boldx , \: \boldy$, where
$\boldy$ is an $n \times 1$ vector of new indeterminates. \\

\noindent 4. Impose SOS constraints on $p(\boldx,\boldy)$, and
$-r(\boldx,\boldy)$, and solve the associated SOS feasibility
problem. If a solution exists, the SOS solver will find values for
the unknown coefficients, such that the constraints are satisfied. \\

\noindent 5. Use the obtained coefficients $a_i, b_i, c_i$ to
construct the contraction metric $\boldM(\boldx)$ and the
corresponding $\boldR(\boldx)$. \\

\noindent 6. Optionally, for graphical presentation, independent
verification, or if the convex optimization procedure runs into
numerical error, further testing can be done to verify the
validity of the computed solution. To do this, we can check if the
matrix constraints $\boldM(\boldx) \succ 0$, and $\boldR(\boldx)
\prec 0$ hold over a range of the state space by finding and
plotting the eigenvalues over this range.  If a true feasible
solution does not exist, the minimum eigenvalue of
$\boldM(\boldx)$ will be negative or the maximum eigenvalue of
$\boldR(\boldx)$ will be positive. Either one of these cases
violates the matrix constraints which certify contraction. In most
semidefinite programming solvers, the matrix $Q$ in (\ref{Q}) is
computed with floating point arithmetic.  If $Q$ is near the
boundary of the set of positive semidefinite matrices, it is
possible for the sign of eigenvalues that are zero or close to
zero to be computed incorrectly from numerical roundoff and for
the semidefinite program solver to encounter numerical
difficulties. Numerical issues are further discussed in Section
\ref{NumericalProblems}. \\

\noindent 7. An explicit lower bound on the rate of convergence
can be found by using bisection to compute the largest $\beta$ for
which there exist matrices $\boldM(\boldx) \succ 0$ and $
\boldR_{\beta}(\boldx) =
\pardfdxnm ^T \boldM + \boldM
\pardfdxnm + \dot{\boldM} + \beta \boldM \prec 0$.\\

For the specific examples presented later in the paper, we have
used SOSTOOLS, a SOS toolbox for MATLAB developed for the
specification and solution of sums of squares programs
\cite{SOSTOOLS}. The specific structure of SOS matrices, or
equivalently, the bipartite form of the polynomials $p(\boldx,
\boldy)$ and $r(\boldx, \boldy)$ is exploited through the option
\texttt{sparsemultipartite} of the command \texttt{sosineq} that
defines the SOS inequalities. Future versions of SOSTOOLS will
allow for the direct specification of matrix SOS constraints.

We present next are two examples of using this procedure to search
for contraction metrics for nonlinear systems with polynomial
dynamics. The systems studied are a model of a jet engine with
controller, and a Van der Pol oscillator.

\subsection{Example: Moore-Greitzer Jet Engine Model}\label{jetexamples}

The algorithm described was tested on the following dynamics,
corresponding to a Moore-Greitzer model of a jet engine, with
stabilizing feedback operating in the no-stall mode \cite{KKK}. In
this model, the origin is translated to a desired no-stall
equilibrium.  The state variables correspond to $\phi = \Phi - 1$,
$\psi = \Psi - \Psi_{co} - 2$, where $\Phi$ is the mass flow, $\Psi$
is the pressure rise and $\Psi_{co}$ is a constant \cite{KKK}.
The dynamic equations take the form:
\begin{equation}\label{JetDynamics}
\left[%
\begin{array}{c}
  \dot{\phi} \\
  \dot{\psi} \\
\end{array}%
\right] =\left[%
\begin{array}{c}
  -\psi - \frac{3}{2}\phi^2 - \frac{1}{2}\phi^3 \\
  3\phi - \psi \\
\end{array}%
\right]
\end{equation}
The only real-valued equilibrium of the system is $\phi = 0$, $ \psi=
0$. This equilibrium is stable.

The results of the algorithmic search for SOS matrices
$\boldM(\boldx)$ and $-\boldR(\boldx)$ of various orders are given
in Table~\ref{JetTable}. Values in the table, except the final
row, are output values from SeDuMi \cite{SeDuMi}, the semidefinite
program solver used as the optimization engine in solving the SOS
program. \texttt{CPU time} is the number of seconds it took for
SeDuMi's interior point algorithm to find a solution. As expected,
the computation time increases with the degree of the polynomial
entries of $\boldM(\boldx)$. \texttt{Feasibility ratio} is the
final value of the feasibility indicator. This indicator converges
to $1$ for problems with a complementary solution, and to $-1$ for
strongly infeasible problems. If the feasibility ratio is
somewhere in between, this is usually an indication of numerical
problems.  The values \texttt{pinf} and \texttt{dinf} detect the
feasibility of the problem.  If $\texttt{pinf}=1$, then the primal
problem is infeasible. If $\texttt{dinf}=1$, the dual problem is
infeasible.  If \texttt{numerr} is positive, the optimization
algorithm (i.e., the semidefinite program solver) terminated
without achieving the desired accuracy. The value
$\texttt{numerr}=1$ gives a warning of numerical problems, while
$\texttt{numerr}=2$ indicates a complete failure due to numerical
problems.

As shown in Table~\ref{JetTable}, for this system no contraction
metric with polynomial entries of degree 0 or 2 could be found.
This can be certified from the solution of the dual optimization
problem. Since SeDuMi is a primal-dual solver, this infeasibility
certificates are computed as a byproduct of the search for
contraction metrics.

An explicit lower bound for the rate of convergence of the
trajectories of the jet engine model, i.e., the largest value
$\beta$ for which matrices $\boldM(\boldx) \succ 0$ and $
\boldR_{\beta}(\boldx) =
\pardfdxnm ^T \boldM + \boldM
\pardfdxnm + \dot{\boldM} + \beta \boldM \prec 0$ were found, was $\beta = 0.818$.

We remark that for this system, it is also possible to prove
stability using standard Lyapunov analysis techniques. However, we
illustrate stability of this example from a contraction viewpoint
because contraction theory offers a good approach to study this
system when there is parametric uncertainty in the plant dynamics
or feedback equations. For example, in the no-stall mode, the jet
dynamics equations are
\begin{equation}\label{JetDynamicscontrol}
\left[%
\begin{array}{c}
  \dot{\phi} \\
  \dot{\psi} \\
\end{array}%
\right] =\left[%
\begin{array}{c}
  -\psi - \frac{3}{2}\phi^2 - \frac{1}{2}\phi^3 \\
  -u\\
\end{array}%
\right]
\end{equation} where $u$ is a control variable.
If a nominal stabilizing feedback control $u$ can be found (e.g.,
using backstepping \cite{KKK} or some other design method), the SOS
techniques described in Section \ref{Uncert-alg-des} provide a way to
find other stabilizing feedback controls which are centered around the
nominal control.  For example, if a stabilizing linear feedback
control $u = k_{1} \phi + k_{2} \psi$ can be found, we can interpret
$k_{1}$ and $k_{2}$ as uncertain parameters and use the methods
described in Section \ref{Uncert-alg-des} to search for ranges of gain
values centered around the nominal values $k_{1}$ and $k_{2}$ that
will also stabilize the system.

\begin{table}
\begin{footnotesize}\begin{center}
\begin{tabular}{|l||c|c|c|c|}
  \hline
  Degree of polynomials in $\boldM(\boldx)$& 0 &2 & 4  & 6  \\
  \hline
  CPU time (sec)  &0.140 & 0.230 & 0.481& 0.671  \\
  \hline
  Feasibility ratio& -1.000  &-0.979 & 1.003 & 0.990  \\
  \hline
  \texttt{pinf} &1 & 1& 0 & 0  \\
  \hline
  \texttt{dinf} & 0& 0 & 0 & 0  \\
  \hline
  \texttt{numerr} &0 & 1 & 0 & 0  \\
  \hline
  $M \succ 0$, $R \prec 0$ conditions met? & no & no & yes & yes \\
  \hline
\end{tabular} \end{center}\end{footnotesize}
\caption{Contraction matrix search results for closed-loop jet
engine dynamics.}\label{JetTable}
\end{table}

\subsection{Example: Van der Pol Oscillator}

A classic example that has played a central role in the
development of nonlinear dynamics is given by the Van der Pol
equation
\begin{equation}\label{vdposc}
\ddot{x} + \alpha(x^2 + k)\dot{x} + \omega^2x = 0,
\end{equation}
with $\alpha \geq 0$, $k$, and $\omega$ as parameters.  Historically
this equation arose from studying nonlinear electric circuits used in
the first radios \cite{Strogatz}. When $k<0$, the solutions of
(\ref{vdposc}) behave like a harmonic oscillator with a nonlinear
damping term $\alpha (x^2 + k)\dot{x}$.  The term provides positive
damping when $|x|>k$ and negative damping when $|x|<k$. Thus, large
amplitude oscillations will decay, but if they become too small they
will grow larger again \cite{Strogatz}. If $k >0$ all trajectories
converge to the origin.

In Table~\ref{OscTableu0} we present the results of running the
contraction matrix search algorithm for the system
\begin{displaymath}
\left[
\begin{array}{c}
  \dot{x}_1 \\
  \dot{x}_2 \\
\end{array}\right] = \left[
\begin{array}{c}
  x_2 \\
  -\alpha (x_1^2 + k)x_2 - \omega^2 x_1  \\
\end{array}
\right], \end{displaymath} with $\alpha = 1, \omega = 1$, which is the
state-space version of the Van der Pol oscillator (\ref{vdposc}). We
present solution for various values of $k$, with a contraction matrix
with entries that are quartic polynomials.

As a natural first step we searched for a constant contraction
metric.   None could be found algorithmically.  This was expected
as it is easily shown analytically that a constant contraction
matrix for this system does not exist. If $\boldM$ is constant,
then {\small
\begin{eqnarray}
\boldM  = \left[%
\begin{array}{cc}
  a & b \\
  b & c \\
\end{array}%
\right], \qquad \pardfdxnm = \left[%
\begin{array}{cc}
  0 & 1 \\
  -1 - 2x_1x_2 & -x_1^2 - k \\
\end{array}%
\right], \nonumber \\ \nonumber \\
\boldR = \left[%
\begin{array}{cc}
  -2b - 4bx_1x_2 & a - bx_1^2 - kb - c - 2cx_1x_2 \\
  a - bx_1^2 - kb - c - 2cx_1x_2 & 2b-2cx_1^2 - 2kc \\
\end{array}%
\right].
\end{eqnarray}} For $\boldR$ to be negative definite $\boldR_{11}$ must be
negative for all values of $x_1$, $x_2$. In other words $-2b -
4bx_1x_2 \leq 0$ or $-1 \leq 2 x_1x_2$.  This clearly does not hold
for all values of $x_1$, and $x_2$. A more complicated analysis (or a
duality argument) also shows why there is no contraction matrix with
quadratic entries for this system.

\begin{table*}
\begin{footnotesize} \begin{center}

\begin{tabular}{|c||c|c|c|c|c||c|c|c|c|c|c|c|}
\hline
  Degree of polynomials in $\boldM(\boldx)$ &4 &4 &4 &4 &4 &4 &4 &4
  &4 &4\\
  \hline
  k   & -10 & -1 & -0.1 & -0.01 & -0.001 & 0.001 & 0.01 & 0.1 & 1 & 10 \\
  \hline
  pinf  & 1 & 1 & 1 & 0 & 0 & 0 & 0 & 0 & 0 & 0 \\
  \hline
  dinf  & 0 & 0 & 0 & 0 & 0 & 0 & 0 & 0 & 0& 0\\
  \hline
  numerr  & 1 & 1 & 1 & 1 & 1 & 0 & 0 & 0 & 0 & 0 \\
  \hline
 $M \succ 0$, $R \prec 0$ conditions met? & no & no & no & no & no & yes & yes & yes & yes & yes \\
  \hline

\end{tabular}
\end{center} \end{footnotesize}
\caption{Contraction matrix search results for oscillator
dynamics.}\label{OscTableu0}
\end{table*}

\begin{figure}
  \begin{center}
    \subfigure[$ k  = 0.5$]{\label{pp_osc_kpos}\includegraphics[scale=0.4]{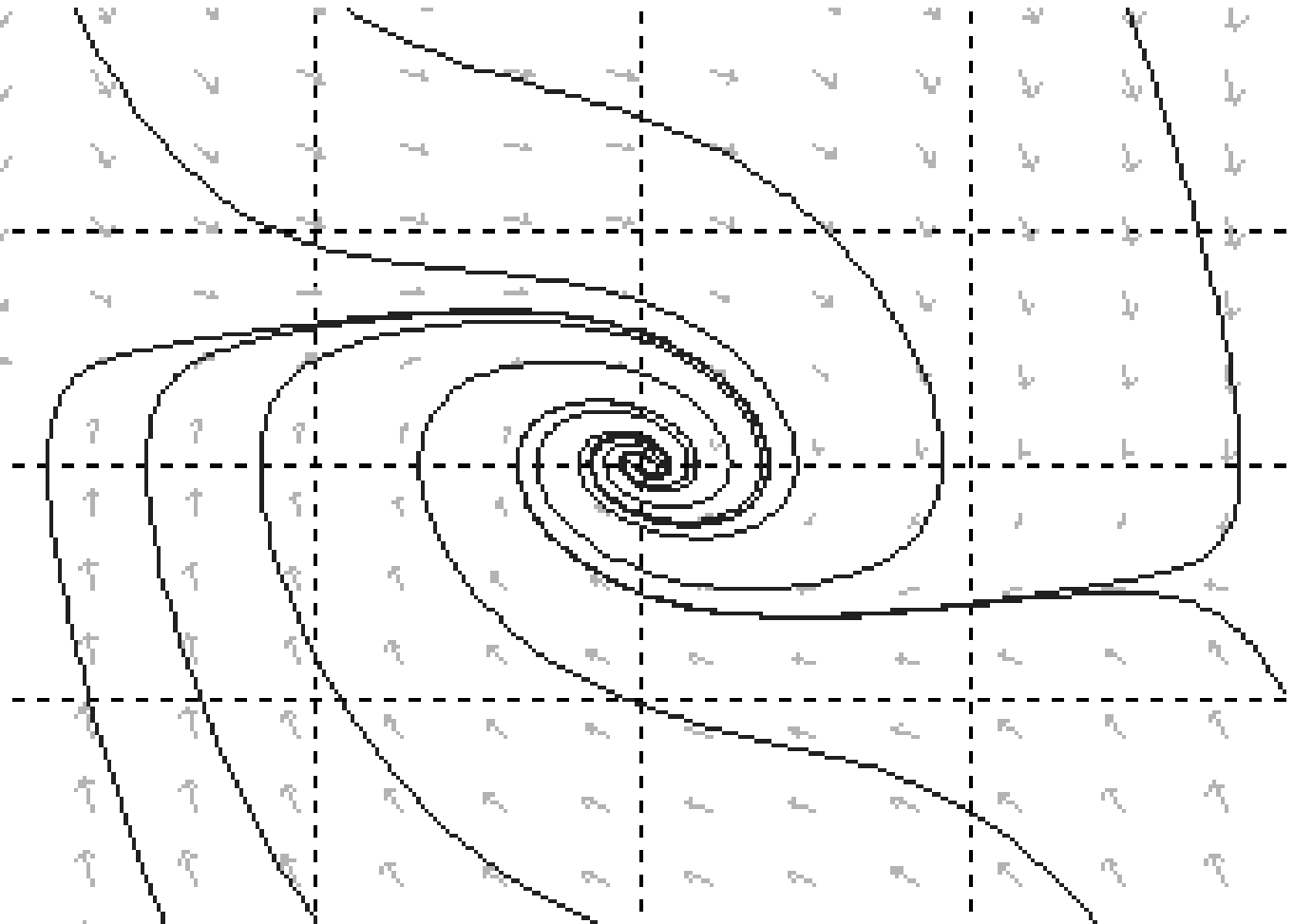}}
   \hspace{1.5 cm}
    \subfigure[$k = -0.5$]{\label{pp_osc_kneg}\includegraphics[scale=0.4]{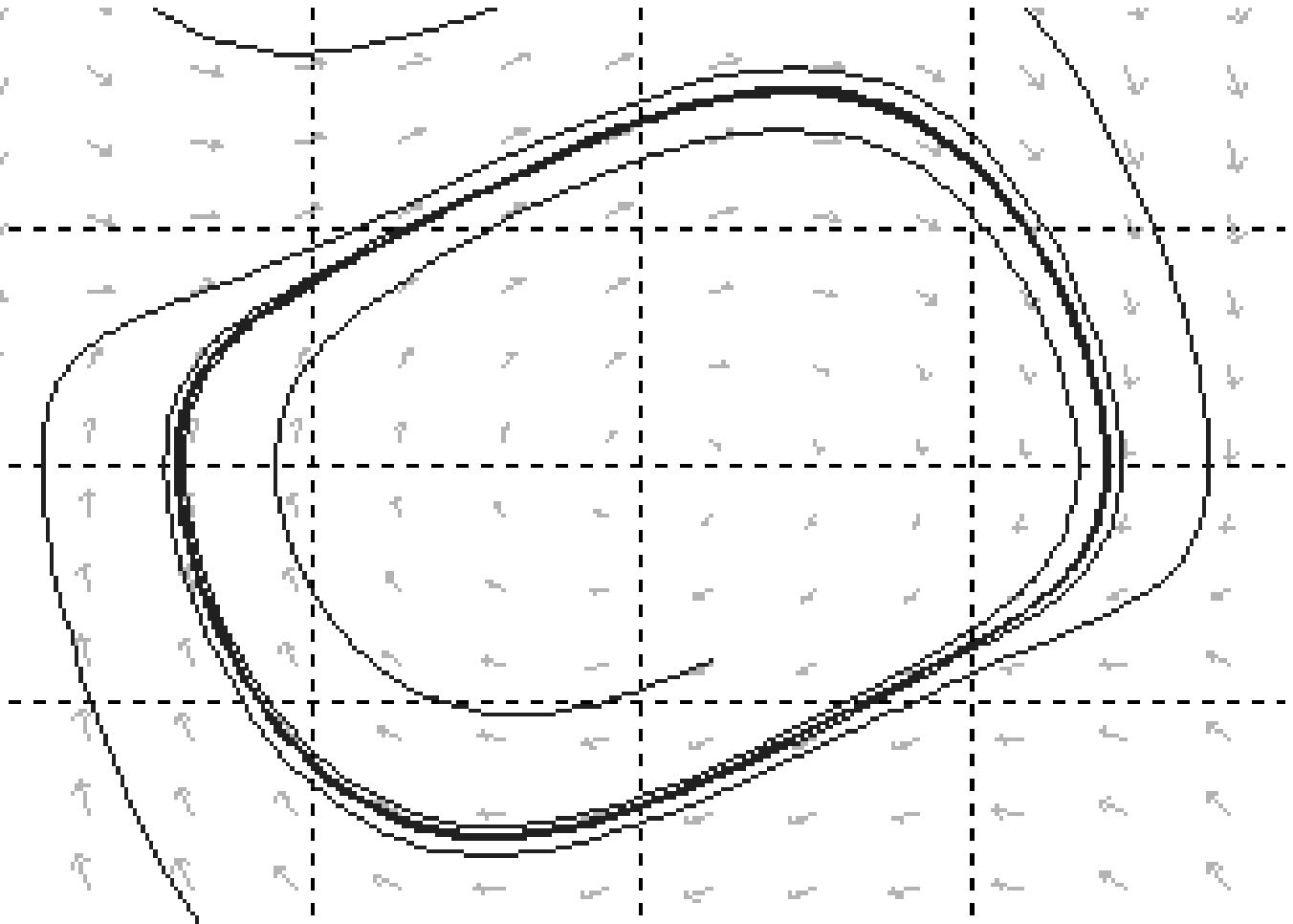}}
  \end{center}
  \caption{Phase plots of Van der pol Oscillator.}
  \label{pp_vdpo}
\end{figure}

The algorithm finds a contraction function for the system
$\ddot{x} + (x^2 + k)\dot{x} + x = 0$ when $k >0$ but not when
$k<0$. As shown in Figure \ref{pp_vdpo} the trajectories of the
oscillator converge to zero when $k > 0$, and converge to a
limit-cycle when $k <0$. Thus, the results of the contraction
metric search is as expected.  Since all trajectories converge to
the origin when $k >0$ we expect that a contraction metric exists
for the system.  In the case where $k <0$ the origin is an
unstable fixed point and thus the system is not contracting.

Since for $k<0$ the system is not contracting, we should not be able
to find a contraction function. It should be noted that the converse
does not hold. The fact that we cannot find a contraction function
does not necessarily mean that the system is not contracting. This is
because finding an SOS representation of the constrained quadratic
functions is a sufficient condition for their positivity, not a
necessary one.

It should be noted that for the example above, we can prove
stability through Lyapunov analysis, and SOS programming can also
be used to find Lyapunov functions \cite{SOSlyap}.  However, we
illustrate this example here as contraction theory applied to a
slightly modified version of this system provides a nice way to
prove synchronization of coupled Van der Pol oscillators. This
will be discussed in Section \ref{EI}. This synchronization
property is much more difficult to prove with standard Lyapunov
methods.

\section{Contraction Metrics and Systems with Uncertain Dynamics}
\label{UD}
\subsection{Uncertainty Analysis with Contraction Metrics and
SOS Programming}\label{Uncert-alg-des}

From the robust control perspective, one of the most appealing
features of contraction theory is the fact that it provides a natural
framework in which to study uncertain nonlinear systems where the
parametric uncertainty changes the location of the equilibrium
points. In general, standard Lyapunov analysis does not handle this
situation particularly well, since the Lyapunov function must track
the changes in the location of the steady-state solutions, thus
forcing the use of parameter-dependent Lyapunov functions.  However,
in general it may be impossible to obtain any kind of closed form
expression of the equilibria in terms of the parameters, thus
complicating the direct parametrization of possible Lyapunov
functions.

Much attention has been given to robust stability analysis of linear
systems (e.g.,
\cite{PDL:Haddad,PDL:Gahinet,PDL:Feron,LMIBook,Zhou}). Less attention,
however, has been paid to nonlinear systems with moving
equilibria. Two papers addressing this issue are
\cite{Wang-Michel,Rantzerpaper}. The approach in \cite{Wang-Michel} is
to consider systems described by the equations
\begin{equation}
\dot{\mathbf{x}} = \mathbf{f}(\mathbf{x}) + \mathbf{h}(\mathbf{x}),
\label{wangeq}
\end{equation}
where $\mathbf{x}$ is a real $n$-vector, $\mathbf{f}$ and $\mathbf{h}$
are continuously differentiable functions, and
$\mathbf{h}(\mathbf{x})$ represents the uncertainties or perturbation
terms. Given an exponentially stable equilibrium $\mathbf{x}_e$,
\cite{Wang-Michel} establishes sufficient conditions by using the
linearization of the system to produce Lyapunov functions which prove
existence and local exponential stability of an equilibrium
$\tilde{\mathbf{x}}_e$ for (\ref{wangeq}) with the property
$|\mathbf{x}_e - \tilde{\mathbf{x}}_e| < \varepsilon$ where
$\varepsilon$ is sufficiently small.

Since the approach in \cite{Wang-Michel} is essentially based on a
fixed Lyapunov function, it is more limited than our approach using
contraction theory and SOS programming, and can prove stability only
under quite conservative ranges of allowable uncertainty.  By
``allowable'' we mean that if the uncertainty is in this range, the
equilibrium remains exponentially stable under the uncertainty. A
quantitative measure of this conservativeness will be given in Section
\ref{Scalar Additive Uncertainty} where we discuss the results of the
method applied to an uncertain model of a Moore-Greitzer jet engine
and compare them to an approach via contraction theory and SOS
programming.

The approach in \cite{Rantzerpaper} is to linearize the dynamics
around an equilibrium which is a function of the uncertain
parameter ($\boldx_0 = \boldg(\delta), \; \delta \in \Omega$) and
then use structured singular values to determine the eigenvalues
of the linearized system $\frac{d\boldz}{dt} = A(\delta)\boldz$ if
$A(\delta)$ is rational in $\delta$.
If $A(\delta)$ is marginally stable, no conclusions can be made
about the stability of the nonlinear system.

The contraction theory framework eliminates the need for
linearization, and even the need to know the exact position of the
equilibrium, in order to analyze stability robustness in uncertain
nonlinear systems. In contrast to the Lyapunov situation, when
certain classes of parametric uncertainty are added to the system,
a contraction metric for the nominal system will often remain a
contraction metric for the system with uncertainty, even if the
perturbation has changed the equilibrium of the nonlinear system.

As noted in Section \ref{CT}, if a global time-invariant
contraction metric exists for an autonomous system, and an
equilibrium point exists for the system, all trajectories converge
to a unique equilibrium point, and we can always produce a
Lyapunov function of the form $V(\boldx) =
\boldf(\boldx)^T\boldM(\boldx)\boldf(\boldx)$. When a system
contains parametric uncertainty, this formula yields the
parameter-dependent Lyapunov function $V(\boldx, \delta) =
\boldf(\boldx, \delta)^T\boldM(\boldx)\boldf(\boldx, \delta)$ for
ranges of the parametric uncertainty $\delta$ where the
contraction metric for the nominal system is still a contraction
metric for the system with perturbed dynamics. Thus, if a
contraction metric can be found for the system, we can easily
construct a Lyapunov function which tracks the uncertainty for a
certain range.

\subsubsection{Case 1: Bounds on the uncertainty range
for which the system remains contractive with respect to the nominal
metric.}

We can estimate the range of uncertainty under which the
contraction metric for the nominal system is still a contraction
metric for the perturbed system. To calculate this range, a SOS
program can be written to minimize or maximize the amount of
uncertainty allowed subject to the constraint
$\boldR_{\delta}(\boldx) =
\pardfdxdel ^T \boldM + \boldM
\pardfdxdel + \dot{\boldM}(\boldf_{\delta}(\boldx)) \prec 0$, where $\boldf_{\delta}(\boldx)$ are the dynamics
for the system with parametric uncertainty.  The uncertainty bound
is a decision variable in the SOS program and enters the
constraint above in the $\pardfdxdel$ and
$\boldf_{\delta}(\boldx)$ terms.

If we have more than one uncertain parameter in the system, we can
find a polytopic inner approximation of the set of allowable
uncertainties with SOS Programming. For example, if we have two
uncertain parameters, we can algorithmically find a polytope in
parameter space for which the original metric is still a contraction
metric. The convex hull of four points, each which can be found by
entering one of the four combinations, $(\delta _1, \delta_2) =
(\gamma, \gamma)$, $(\delta _1, \delta_2) = (\gamma, -\gamma)$,
$(\delta _1, \delta_2) = (-\gamma, \gamma)$, or $(\delta _1, \delta_2)
= (-\gamma, -\gamma)$, into the uncertainty values in
$\boldf_{\mathbf{\delta} = [\delta_1, \delta_2]^T}(\boldx)$ and then
maximizing $\gamma$ subject to the constraint $\boldR_{\gamma}(\boldx)
= \pardfdxgam ^T \boldM + \boldM \pardfdxgam +
\dot{\boldM}(\boldf_{\gamma}(\boldx)) \prec 0$, defines a polytope
over which stability is guaranteed.

\subsubsection{Case 2: Search for a contraction metric that guarantees the largest
symmetric uncertainty interval for which the system is contractive.}

Alternatively, we can instead optimize the search for a metric
$\boldM(\boldx)$ that provides the largest symmetric uncertainty
interval for which we can prove the system is contracting. If the
scalar uncertainty $\delta$ enters the system dynamics affinely, in
other words if $\boldf(\boldx) = \mathbf{f_1}(\boldx) + \delta
\, \mathbf{f_2}(\boldx)$, we can perform this optimization as
follows. First write $\boldR(\boldx, \delta) = \mathbf{R_0}(\boldx) +
\delta \mathbf{R_1}(\boldx)$.  To find the largest interval $(-\gamma,
\gamma)$ such that for all $\delta$ that satisfy $-\gamma < \delta <
\gamma$ the system is contracting, introduce the following constraints
into an SOS program: \begin{displaymath} \boldM(\boldx) \succ 0, \quad
\mathbf{R_0}(\boldx) + \gamma \mathbf{R_1}(\boldx) \prec 0, \quad
\mathbf{R_0}(\boldx) -\gamma \mathbf{R_1}(\boldx) \prec 0.
\end{displaymath}  We note that $\gamma$ multiplies the scalar decision coefficients
$a_i$, $b_i$, and $c_i$ in $\boldR_1(\boldx)$ and thus we must use
a bisection procedure to find the maximum value of $\gamma$ for
which there exists SOS matrices $\boldM(\boldx)$,
$\mathbf{R_0}(\boldx)$ and $\mathbf{R_1}(\boldx)$ that satisfy the
constraints above.

If there are two uncertain parameters that enter the system
dynamics affinely, we can extend the procedure above as follows:
To find the largest uncertainty square with width and height
$\gamma$ such that for all $\delta_1$ and $\delta_2$ that satisfy
$-\gamma < \delta_1 < \gamma$ and $-\gamma < \delta_2 < \gamma$
the system is contracting, first write $\boldR(\boldx, \delta_1,
\delta_2) = \mathbf{R_0}(\boldx) + \delta_1 \mathbf{R_1}(\boldx) +
\delta_2 \mathbf{R_2}(\boldx)$,  Then introduce the following
constraints into and SOS program:
\begin{eqnarray}
\boldM(\boldx) \succ 0, \quad \mathbf{R_0}(\boldx) + \gamma
\mathbf{R_1}(\boldx) + \gamma \mathbf{R_2}(\boldx) \prec 0, \quad
\mathbf{R_0}(\boldx) + \gamma \mathbf{R_1}(\boldx) - \gamma \mathbf{R_2}(\boldx) \prec 0 \nonumber \\
\mathbf{R_0}(\boldx) - \gamma \mathbf{R_1}(\boldx) + \gamma
\mathbf{R_2}(\boldx) \prec 0, \quad \mathbf{R_0}(\boldx) - \gamma
\mathbf{R_1}(\boldx) - \gamma \mathbf{R_2}(\boldx) \prec 0.
\end{eqnarray}
Next, as in the scalar uncertainty case, use a bisection procedure to
find the maximum value of $\gamma$ for which there exists SOS matrices
$\boldM(\boldx)$, $\mathbf{R_0}(\boldx) $, $\mathbf{R_1}(\boldx)$ and
$\mathbf{R_2}(\boldx)$ that satisfy the constraints above.  In the
case of a large number of uncertain parameters, standard relaxation
and robust control techniques can be used to avoid an exponential
number of constraints.

\subsection{Example: Moore-Greitzer Jet Engine Model with Uncertainty}
\subsubsection{Scalar Additive Uncertainty}\label{Scalar Additive
Uncertainty}

As described above, SOS programming can be used to find ranges of
uncertainty under which a system with uncertain perturbations is
still contracting with the original contraction metric. The
contraction metric found for the deterministic system continues to be
a metric for the perturbed system over a range of uncertainty even if
the uncertainty shifts the equilibrium point and trajectories of the
system. For the Moore-Greitzer jet engine model, the dynamics in
(\ref{JetDynamics}) were perturbed by adding a constant term $\delta$
to the first equation.
\begin{equation}\label{JetDynamics_eps}
\left[%
\begin{array}{c}
  \dot{\phi} \\
  \dot{\psi} \\
\end{array}%
\right] =\left[%
\begin{array}{c}
  -\psi - \frac{3}{2}\phi^2 - \frac{1}{2}\phi^3 + \delta \\
  3\phi - \psi \\
\end{array}%
\right]
\end{equation}

\begin{table}
\begin{center}
\footnotesize
\begin{tabular}[h]{|l||c|c|c|c|}
  \hline
  Degree of polynomials in $\boldM(\boldx)$ & 4  & 6  \\
  \hline
  $ \delta$ range &(-0.126,0.630)  &( -0.070, 0.635) \\
    \hline
\end{tabular} \normalsize \end{center}
\caption{Range of perturbation where closed-loop uncertain jet
engine dynamics given in (\ref{JetDynamics_eps}) are contracting
with respect to the nominal metric.} \label{epsJettable}
\end{table}

In Table \ref{epsJettable} we display the ranges of $\delta$ where
the system was still contracting with the original contraction
metric for $4^{th}$ and $6^{th}$ degree contraction metrics. Note
the range of allowable uncertainty is not symmetric.

\begin{figure}[]
  \begin{center}
    \subfigure[$\delta =
    -0.5$]{\label{JetDelMpt5}\includegraphics[scale=0.5]{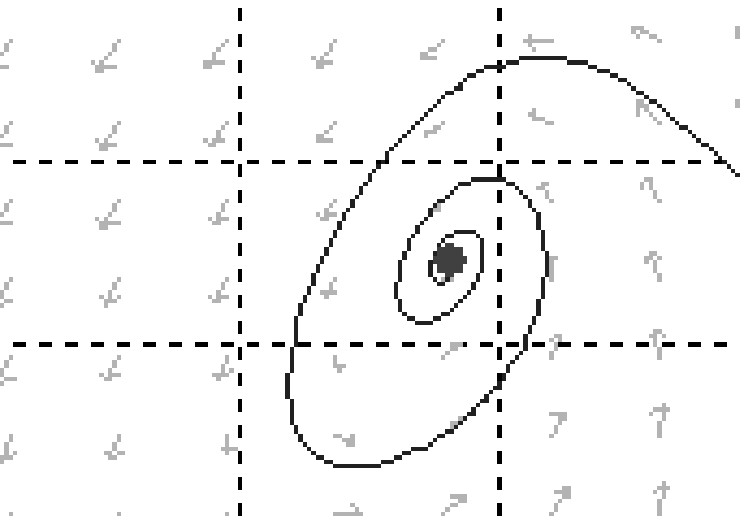}}
    \hspace{1cm}
    \subfigure[$\delta = -1.01$]{\label{JetDelM1pt01}\includegraphics[scale=0.5]{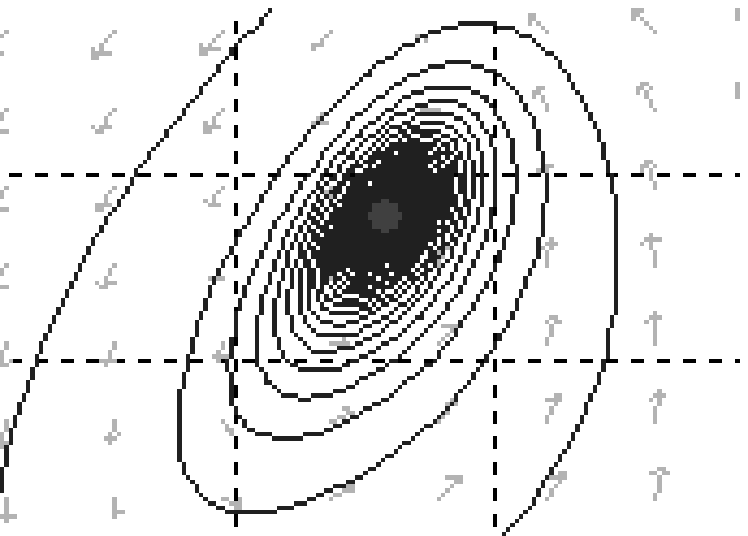}}
    \hspace{1cm}
    \subfigure[$\delta = -1.1$]{\label{JetDelM1pt1}\includegraphics[scale =0.5]{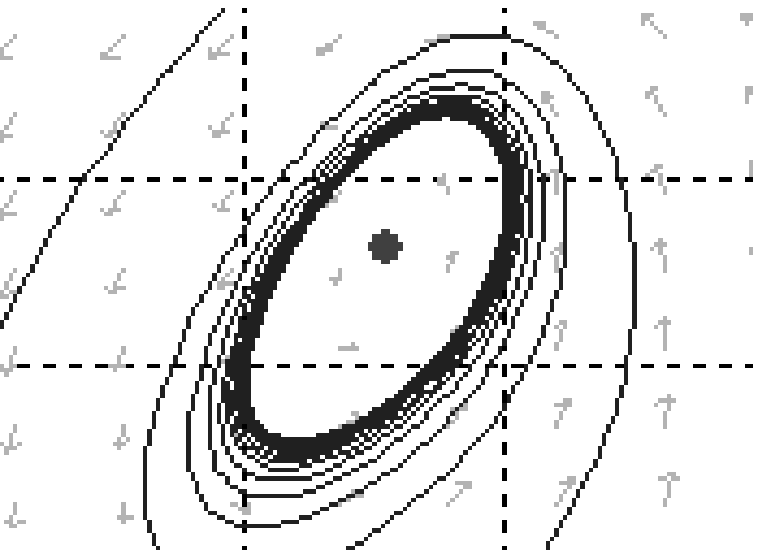}}
  \end{center}
  \caption{Hopf bifurcation in uncertain jet dynamics.}
  \label{fig:Hopf}
\end{figure}

When instead we optimized the contraction metric search to get the
largest symmetric $\delta$ interval we obtained the results listed
in Table~\ref{tbl:optdelta}.  A $6^{th}$ degree contraction
function finds the uncertainty range $|\delta| \leq 1.023$.
Because a Hopf bifurcation occurs in this system at $\delta
\approx 1.023$, making the system unstable for $\delta > 1.023$,
we can conclude that the 6th degree contraction metric is the
highest degree necessary to find the maximum range of uncertainty
for which the system is contracting. The Hopf bifurcation is shown
in Figure~\ref{fig:Hopf}.

Using the techniques in \cite{Wang-Michel} we computed the
allowable uncertainty range for the system given in
(\ref{JetDynamics_eps}) as $|\delta| \leq 5.1 \times 10^{-3}$. In
the notation of \cite{Wang-Michel}, we calculated the other
parameters in Assumption 1  of \cite{Wang-Michel} as: $\mathbf{h}
= [\delta, \; \; 0]^T$, $|A^{-1}|_\infty = 1$, $|D
\mathbf{h}(\boldx_e)|_\infty = 0$, $a = \frac{1}{30}$, and
$|\mathbf{h}(\boldx_e)|_\infty = |\delta|$, where $\delta$ is the
perturbation term in (\ref{JetDynamics_eps}).  The allowable range
$|\delta| \leq 1.023$ computed via contraction theory and SOS
programming is much larger than the allowable uncertainty range
$|\delta| \leq 5.1 \times 10^{-3}$ computed with the techniques in
\cite{Wang-Michel}.

\begin{table*}
\begin{footnotesize}\begin{center}
\begin{tabular}{|l||c|c|c|}
  \hline
  Degree of polynomials in $\boldM(\boldx)$ & 4  & 6 & 8  \\
  \hline
   $ \delta$ range  &$|\delta| \leq 0.938$  & $|\delta| \leq 1.023$ & $|\delta| \leq 1.023$ \\
  \hline
\end{tabular} \end{center}\end{footnotesize}
\caption{Symmetric range of perturbation for the jet engine
model.} \label{tbl:optdelta}
\end{table*}

\subsubsection{Scalar Multiplicative Uncertainty}

The approaches in Section \ref{Uncert-alg-des} also apply to
multiplicative uncertainty, since the multiplicative coefficients
enter affinely in the constraints in the SOS program. Tables
\ref{epsJettable-mult} and \ref{tbl:optdelta:mult} present the results
of the described uncertainty analysis on the following system, which
is equation (\ref{JetDynamics}) with multiplicative uncertainty.
\begin{equation}\label{JetDynamics:multuncert}
\left[%
\begin{array}{c}
  \dot{\phi} \\
  \dot{\psi} \\
\end{array}%
\right] =\left[%
\begin{array}{c}
  -\psi - \frac{3}{2}\phi^2 - \frac{1}{2} \delta \phi^3 \\
  3\phi - \psi \\
\end{array}%
\right].
\end{equation}

\begin{table}
\begin{center}
\footnotesize
\begin{tabular}{|l||c|c|c|c|}
  \hline
  Degree of polynomials in $\boldM(\boldx)$ & 4  & 6  \\
  \hline
  $ \delta$ range & (0.9767, 5.8686) & (0.9796, 3.9738)  \\
    \hline
\end{tabular} \normalsize \end{center}
\caption{Range of perturbation for which the uncertain system given in
(\ref{JetDynamics:multuncert}) is contracting with respect to the
nominal metric.}
\label{epsJettable-mult}
\end{table}

\begin{table*}
\begin{footnotesize}\begin{center}
\begin{tabular}{|l||c|c|c|}
  \hline
  Degree of polynomials in $\boldM(\boldx)$ & 4  & 6 & 8  \\
  \hline
   $ \delta$ range  &$(1- 0.247,1+ 0.247)$  & $(1 - 0.356  ,1+ 0.356 )$ & $(1 - 0.364  ,1+0.364)$ \\
  \hline
\end{tabular} \end{center}\end{footnotesize}
\caption{Symmetric range of perturbation where uncertain
closed-loop jet engine dynamics given in
(\ref{JetDynamics:multuncert}) are contracting.}
\label{tbl:optdelta:mult}
\end{table*}

\subsubsection{Multiple Uncertainties}

\begin{figure}[]
  \begin{center}
   \includegraphics[scale=0.4]{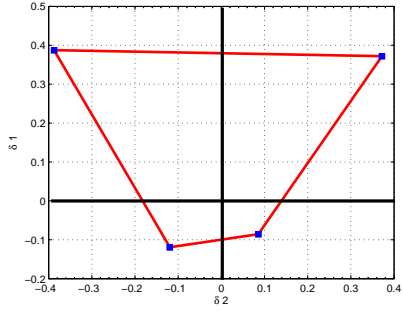}
   \end{center}
  \caption{Polytopic region of uncertainty where closed-loop jet engine
dynamics given in (\ref{JetDynamics:2Duncert}) are contracting with respect to nominal metric.}
  \label{poly-uncert}
\end{figure}

We consider next the system that results from introducing two additive
uncertainties to the jet dynamics in equation (\ref{JetDynamics}). We
computed an uncertainty polytope (shown in Figure \ref{poly-uncert})
for which the system
\begin{equation}\label{JetDynamics:2Duncert}
\left[%
\begin{array}{c}
  \dot{\phi} \\
  \dot{\psi} \\
\end{array}%
\right] =\left[%
\begin{array}{c}
  -\psi - \frac{3}{2}\phi^2 - \frac{1}{2}\phi^3 + \delta_1 \\
  3\phi - \psi + \delta_2 \\
\end{array}%
\right]
\end{equation}
is guaranteed to be contracting with respect to the original metric.
 Table \ref{tbl:optdelta:2D} shows the results of optimizing the
contraction metric to find the largest uncertainty square with
width and height $\gamma$ such that for all $\delta_1$ and
$\delta_2$ that satisfy $-\gamma < \delta_1 < \gamma$ and $-\gamma
< \delta_2 < \gamma$ the system is contracting.
\begin{table*}
\begin{footnotesize}\begin{center}
\begin{tabular}{|l||c|c|c|}
  \hline
  Degree of polynomials in $\boldM(\boldx)$ & 4  & 6 & 8  \\
  \hline
    Height and width of allowed uncertainty box  & 0.7093  & 0.7321  & 0.7346  \\
  \hline
\end{tabular} \end{center}\end{footnotesize}
\caption{Symmetric range of perturbation where uncertain
closed-loop dynamics given in (\ref{JetDynamics:2Duncert}) are
contracting.} \label{tbl:optdelta:2D}
\end{table*}

\section{Contraction Metrics for Systems with External
Inputs}\label{EI}
\subsection{Stability Analysis of Systems with External Inputs}

Another interesting feature of the contraction framework is the
relative flexibility in incorporating inputs and outputs. For
instance, to prove contraction of a class of systems with external
inputs, it is sufficient to show the existence of a polynomial
contraction metric with a certain structure. This is described in
the following theorem. \\

\begin{theorem}\label{Thm:EI}
Let
\begin{eqnarray}\label{Thm2:eqn}
 \nonumber \dot{x}_1 &=& f_1(x_1, x_2, \ldots, x_n) \\
\nonumber  \vdots &=& \vdots \\
 \nonumber \dot{x}_k &=& f_k(x_1, x_2, \ldots, x_n) \\
 \nonumber \dot{x}_{k+1} &=& f_{k+1}(x_1, x_2, \ldots, x_n) + v_{k+1}(u) \\
 \nonumber \vdots &=& \vdots\\
   \dot{x}_n &=& f_n(x_1, x_2, \ldots, x_n)+ v_n(u)
\end{eqnarray}
be a set of nonlinear coupled differential equations where only the
last $n-k$ depend explicitly on $u(t)$.  If there exists a $n \times
n$ matrix $\boldM(x_1, \ldots , x_k)$ such that $\boldM \succ 0$ and
$\dot{\boldM} + \frac{\partial f}{\partial x}^T\boldM + \boldM
\frac{\partial f}{\partial x} \prec 0$ then the system is contracting
for all possible choices of $u(t)$.
\label{thm:inputs}
\end{theorem}

\begin{proof}
For notational convenience, let $\dot{\mathbf{x}}_1 = [
\begin{array}{ccc}
  \dot{x}_1 & ... & \dot{x}_k \\
\end{array}]^T = \mathbf{f_1}(\mathbf{x_1}, \mathbf{x_2})$ and
$\dot{\mathbf{x}}_2 = [ \begin{array}{ccc}
  \dot{x}_{k+1} & ... & \dot{x}_n \\
\end{array}]^T = \mathbf{f_2}(\mathbf{x_1}, \mathbf{x_2}, u)$. The metric
$\boldM(x_1, x_2 ..., x_k) = \boldM(\mathbf{x_1})$ is independent
of $\mathbf{x_2}$, and thus $\frac{\partial \boldM_{ij}}{\partial
\mathbf{x_2}} = \mathbf{0} \; \; \forall i,j.$ Since
$\frac{\partial \boldM_{ij}}{\partial t}$ also vanishes, it
follows that $\forall i, j$, $\dot{\boldM}_{ij} = \frac{\partial
\boldM_{ij}}{\partial \mathbf{x_1}}\frac{d \mathbf{x_1}}{dt} +
\frac{\partial \boldM_{ij}}{\partial \mathbf{x_2}}\frac{d
\mathbf{x_2}}{dt} + \frac{\partial \boldM}{\partial t} =
\frac{\partial \boldM_{ij}}{\partial \mathbf{x_1}}\frac{d
\mathbf{x_1}}{dt}$. Thus $\dot{\boldM}(\mathbf{x_1})$ is not a
function of $u(t)$. In addition, $\frac{\partial \boldf}{\partial
\boldx}$ has no dependence on $u(t)$ because $\boldf(\boldx, u) =
\mathbf{h}(\boldx) + \mathbf{v}(u)$. Thus, if there exists a $n
\times n$ matrix $\boldM(x_1, ..., x_k)$ such that $\boldM \succ
0$ and $\dot{\boldM} + \pardfdxnm^T\boldM + \boldM
\pardfdxnm \prec 0$, then the system in (\ref{Thm2:eqn}) is contracting for any value of
$u(t)$.
\end{proof}

Through the example considered in the following section, we will
illustrate how Theorem~\ref{Thm:EI} is particularly useful in proving
synchronization of nonlinear oscillators, an issue explored in more
detail in \cite{OPC-NO}. Theorem~\ref{thm:inputs} can be easily
extended to the case where $u(t)$ is a vector (i.e.  $\boldu(t) =
[u_1(t), ... , u_m(t)]^T$).

\subsection{Coupled Oscillators}
Contraction Theory is a useful tool to study synchronization behaviors
of various configurations of coupled oscillators.  For simplicity, we
only consider here the case a pair of unidirectionally coupled
oscillators; more complicated and general couplings are discussed in
\cite{OPC-NO}.

A state-space model of two unidirectionally coupled oscillators (only
one oscillator influences the other) is
\begin{equation}
\begin{aligned}
    \dot{\boldx} &= \boldf(\boldx,t) \\
    \hfill \dot{\boldy} &= \boldf(\boldy,t) + \boldu(\boldx) -
    \boldu(\boldy),
\end{aligned}
\label{OPC-thm2-eqn}
\end{equation}
where $\boldx, \; \boldy \in \mathbb{R}^m$, are the state vectors,
$\boldf(\boldx,t)$ and $\boldf(\boldy,t)$ are the dynamics of the
uncoupled oscillators, and $\boldu(\boldx) - \boldu(\boldy)$ is
the coupling force\footnote{An example of coupled oscillators
whose state-space representation is in this form is
$\begin{cases}
  \ddot{x} + \alpha(x^2 + k)\dot{x} + \omega^2x  &= 0 \\
 \hfill \ddot{y} + \alpha(y^2 + k)\dot{y} + \omega^2y  &= \alpha \eta (\dot{x} - \dot{y}) \\
\end{cases}$
where $\alpha >0$, $\omega > 0$, $k$  are arbitrary constants.}.
The following theorem is a slightly modified version of Theorem 2
in \cite{OPC-NO}.

\begin{theorem}\label{thm:OPC2}
If $\dot{\boldy} = \boldf(\boldy) + \boldu(\boldy) -
\boldu(\boldx)$ in (\ref{OPC-thm2-eqn}) is contracting with
respect to  $\boldy$ over the entire state space for arbitrary
$\boldu(\boldx)$\footnote{By contracting with respect to $\boldy$
for arbitrary $\boldu(\boldx)$ we mean that the system
$\dot{\boldy} = \boldf(\boldy) - \boldu(\boldy) + \boldu(\boldx)$,
where $\boldy$ is the state vector and $\boldu(\boldx)$ is an
arbitrary driving function, is contracting for all inputs
$\boldu(\boldx)$.}, the two systems will reach synchrony (i.e.
$\boldy(t)$ and $\boldx(t)$ will tend toward the same trajectory)
regardless of initial conditions.
\end{theorem}

\begin{proof}
The system $\dot{\boldy} = \boldf(\boldy) - \boldu(\boldy) +
\boldu(\boldx)$ with input $\boldu(\boldx)$ is contracting with
respect to $\boldy$ over the entire state space and $\boldy (t) =
\boldx (t)$ is a particular solution.  Thus, by the properties of
contraction, all solutions converge exponentially to $\boldy (t) =
\boldx (t)$.
\end{proof}

Theorem \ref{thm:inputs} becomes especially powerful when the
vector field appearing in the second subsystem of
(\ref{OPC-thm2-eqn}) has the structure described in
equation (\ref{Thm2:eqn})\footnote{If it does not have such a structure
and $\boldu(\boldx)$ drives each component of $\boldy$, we lose
degrees of freedom in the possible forms of our contraction
metric.  If $\boldu(\boldx)$ drives each component of $\boldy$ the
only possible contraction metric is a constant.}. We illustrate
this in the next example.

\subsection{Example: Coupled Van der Pol Oscillators}\label{coupled_vdpo}

Consider two identical Van der Pol oscillators coupled as
\begin{equation}\label{coupled-osc-unid}
\begin{cases}
  \ddot{x} + \alpha(x^2 + k)\dot{x} + \omega^2x  &= 0 \\
 \hfill \ddot{y} + \alpha(y^2 + k)\dot{y} + \omega^2y  &= \alpha \eta (\dot{x} - \dot{y}) \\
\end{cases}\end{equation}
where $\alpha >0$, $\omega > 0$, $k$  are arbitrary constants. We
note that if $k<0$, trajectories of the individual oscillator
dynamics converge to a limit cycle. See Figure \ref{pp_osc_kneg}.
We first write these coupled systems in state-space form to get
the equations in the form of (\ref{OPC-thm2-eqn}). Their
state-space form is
\begin{equation}\label{cou-ssform}
\begin{cases}
\left[%
\begin{array}{c}
  \dot{x}_1 \\
  \dot{x}_2 \\
\end{array}%
\right] &= \left[
\begin{array}{c}
  x_2 \\
  -\alpha (x_1^2 + k)x_2 - \omega^2 x_1  \\
\end{array}%
\right] \\
\hfill \\
\hfill \left[%
\begin{array}{c}
  \dot{y}_1 \\
  \dot{y}_2 \\
\end{array}%
\right] &= \left[
\begin{array}{c}
  y_2 \\
  -\alpha (y_1^2 + k + \eta)y_2 - \omega^2 y_1  + \alpha \eta x_2 \\
\end{array}%
\right]. \\
\end{cases}
\end{equation} By Theorem \ref{Thm:EI}, this pair of unidirectional oscillators
will reach synchrony regardless of initial conditions if
\begin{equation}\label{co:eqn4cont}
\dot{\boldy} = \boldf(\boldy) - \boldu(\boldy)+ \boldu(\boldx) = \left[%
\begin{array}{c}
  y_2 \\
  -\alpha(y_1^2 +k + \eta)y_2 - w^2y_1 + \alpha \eta x_2 \\
\end{array}%
\right] \end{equation} is contracting with respect to $\boldy$ for
arbitrary values of $\boldu(\boldx) = x_2$. We see by Theorem
\ref{thm:inputs} that for this to occur, we must find a
contraction metric $\boldM(\boldy)$ that is only  a function of
$y_1$ (i.e. $\boldM(\boldy) = \boldM(y_1)$).

When the search algorithm described in Section~\ref{SearchAlgorithm},
was applied to find a metric that satisfied $\boldM(\boldy) =
\boldM(y_1)$ as well as $\boldM(\boldy)\succ 0$ and $\boldR(\boldy)
\prec 0$, none were found.  However, it is shown in the appendix,
which is a modified version of the appendix of \cite{OPC-NO}, that a
metric that satisfies $\boldM(\boldy)\succ 0$ and $\boldR(\boldy)
\preceq 0$ implies asymptotic convergence of trajectories of system
(\ref{co:eqn4cont}). A system with this metric that satisfies
$\boldM(\boldy)\succ 0$ and $\boldR(\boldy) \preceq 0$ is called
\emph{semi-contracting} \cite{OCA, OPC-NO}.

\noindent The metric
\begin{equation}\label{MAnalyOsc}\boldM(\boldy) =
\left[\begin{array}{cc}
  \omega^2 + \alpha^2(y_1^2 + k + \eta)^2    & \alpha(y_1^2 +k + \eta) \\
  \alpha(y_1^2 + k +\eta ) & 1 \\
\end{array}\right]\end{equation} that appears in \cite{OPC-NO} is only a function of $y_1$
and satisfies $\boldM(\boldy)\succ 0$ and $\boldR(\boldy) \preceq
0$ for the system dynamics (\ref{co:eqn4cont}) if $\alpha >0$ and
$(k + \eta) \geq 0$.  For this $\boldM$ and the system equation
(\ref{co:eqn4cont}), we have
\begin{equation}\label{Rnou} \boldR = \dot{\boldM}  + \pardfdynm^T \boldM + \boldM \pardfdynm = \left[%
\begin{array}{cc}
  -2\alpha \omega^2 y_1^2 - 2\alpha \omega^2(k + \eta ) & 0 \\
  0 & 0 \\
\end{array}%
\right]. \end{equation} For $\alpha >0$, $(k + \eta)>0$,
$\boldM(\boldy) \succ 0$ and $\boldR(\boldy) \preceq 0$. Since
(\ref{MAnalyOsc}) and (\ref{Rnou}) show analytically that the
system (\ref{co:eqn4cont}) is semi-contracting we used our search
algorithm to search for a metric with $\boldM(\boldy)\succ 0$ and
$\boldR(\boldy) \preceq 0$.

\subsubsection{Search for a Semidefinite $\boldR$ Matrix: Numerical
Problems and Solutions}\label{NumericalProblems}

A minor problem that one may encounter when searching for contraction
metrics, depending on the structure of polynomial constraints, is that
the resulting optimization problem may be feasible, but not strictly
feasible. This can cause numerical difficulties in the algorithms used
in the solution procedure. In many cases, however, this can be
remedied by a introducing a presolving stage in which redundant
variables are eliminated. When we ran the search algorithm based on
Theorem~\ref{thm:inputs} and only searched for $\boldM$ as a function
of $y_1$, no valid solution was found even if we only constrained
$\boldR$ to be negative semidefinite and not strictly negative
definite. Since the analytic solution (\ref{Rnou}) was feasible but
not strictly feasible, we hypothesized there was numerical error in
the algorithm. Based on knowledge of the analytic solution
(\ref{Rnou}), we thus constrained $R_{22} = 0$ and $R_{12} = 0$,
eliminated redundant variables, and then searched for a solution in
the resulting lower dimensional space\footnote{Setting $R_{22} = 0$,
and $R_{12} = 0$ leads to redundant decision coefficients in the
polynomial entries of $\boldM$ and $\boldR$.  If these redundant
variables are eliminated through a presolving stage, the search
algorithm finds $\boldM \succ 0$ and $\boldR \preceq 0$.}. With these
constraints in place, a solution was found with the search algorithm.

\section{Conclusions}\label{Conclusions}

In this paper we have described how SOS programming enables an
algorithmic search for contraction metrics for the class of
nonlinear systems with polynomial dynamics.  We also have
illustrated the results through several examples.

These examples illustrate how contraction analysis offers several
significant advantages when compared with traditional Lyapunov
analysis. Contraction analysis provides relative flexibility in
incorporating inputs and outputs. It is also particularly useful
in the analysis of nonlinear systems with uncertain parameters
where the uncertainty changes the equilibrium points of the
system. It is often the case that if the nominal system is
contracting with respect to a metric, the uncertain system with
additive or multiplicative uncertainty will still be contracting
with respect to the original metric, even if the perturbation
changes the equilibrium of the system.  In addition, a slightly
modified version of the standard algorithmic search allows us to
optimize the search to obtain a contraction metric that provides
the largest uncertainty interval for which we can prove the system
is contracting.

Subjects of future research include a careful evaluation of how the
computational resources needed by the algorithm scale with system
size, as well as the benefits and limitations of this approach in the
context of other nonlinear system analysis techniques.

\appendix

\section{Proving Asymptotic Convergence of Coupled Van der Pol Oscillators With a Negative
Semidefinite $\boldR$ Matrix.}

This appendix is a modified version of the appendix in
\cite{OPC-NO}. Consider the system given in (\ref{co:eqn4cont}).
Consider a $2 \times 2$ matrix $\boldM(\boldy)$ that is uniformly
positive definite, and a corresponding $\boldR(\boldy)$ matrix that is
uniformly negative semidefinite, but not uniformly negative
definite. Since (\ref{co:eqn4cont}) is a two-dimensional system, we
can assume without loss of generality that $\boldR(\boldy)$ is of the
form \[\boldR(\boldy) = \left(%
\begin{array}{cc}
  -K(\boldy) & 0 \\
  0 & 0 \\
\end{array}
\right) \]where $K(\boldy)>0 \;\forall \; \boldy$. Let $\delta
\boldy
= \left(%
\begin{array}{cc}
  \delta y_1 & \delta y_2 \\
\end{array}%
\right)^T = \left ( \begin{array}{cc}
  \delta y & \delta \dot{y} \\
\end{array} \right)^T$.  where $y_1$ and $y_2$ are the variables in equation (\ref{co:eqn4cont}). With this $\boldR(\boldy)$ and $\boldM(\boldy)$ matrices, the
general definition of differential length given in
(\ref{RiemannSp}) and associated equation for rate of change of
length (\ref{contractcond2}) are
\[ \delta \boldz^T \delta \boldz = \delta \boldy^T \boldM (\boldy)
\delta \boldy \] and
\begin{eqnarray} \frac{d}{dt} \delta \boldz^T
\delta \boldz &=& \frac{d}{dt}(\delta \boldy^T \boldM (\boldy)
\delta \boldy) \nonumber \\
& = & \delta \boldy ^T \boldR(\boldy) \delta \boldy \nonumber \\
& = & -K(\boldy) \delta y_1^2. \label{negsdR}
\end{eqnarray} Since $\frac{d}{dt}(\delta \boldy^T \boldM (\boldy)
\delta \boldy) \leq 0$ and $\delta \boldz ^T \delta \boldz \geq 0$,
$\delta \boldz ^T \delta \boldz$ has a limit at $t$ goes to
infinity. We will prove through a Taylor series argument that all
trajectories of this system converge asymptotically. If $\delta y =
\delta y_1 \neq 0$, then

\begin{displaymath} \delta \boldz ^T \delta \boldz (t + dt) - \delta \boldz ^T
\delta \boldz (t) = -K(\boldy) (\delta y_1)^2 dt + O((dt)^2)
\end{displaymath} while
if $\delta y_1 = 0$, \[\delta \boldz ^T \delta \boldz (t + dt) -
\delta \boldz ^T \delta \boldz (t) = -2K(\boldy) (\delta
y_2)^2\frac{dt^3}{3!} + O((dt)^4). \] Since $\delta \boldz ^T
\delta \boldz$ converges, $\delta \boldz ^T \delta
\boldz (t + dt) - \delta \boldz ^T \delta \boldz (t)$ approaches
zero asymptotically and hence $\delta y_1$ and $\delta y_2$ or
equivalently $\delta y$ and $\delta \dot{y}$ both tend to zero.
Thus, for any input $\boldu(\boldx)$ all solutions of system
(\ref{co:eqn4cont}) converge asymptotically to a single trajectory
independent of initial conditions, and the unidirectional
oscillators given in (\ref{cou-ssform}) will reach synchrony
asymptotically regardless of initial conditions.

\pagebreak
\nocite{SOS_LMI_Pablo}
\bibliographystyle{plain}
\bibliography{AylwardParriloSlotineJournal}

\begin{thebibliography}{10}

\bibitem{Rantzerpaper}
L.~Andersson and A.~Rantzer.
\newblock Robustness of equilibria in nonlinear systems.
\newblock In {\em Preprints 14th World Congress of IFAC}, volume~E, pages
  129--134, Beijing, P.R. China, 1999.

\bibitem{LMIBook}
S.~Boyd, L.~El Ghaoui, E.~Feron, and V.~Balakrishnan.
\newblock {\em Linear matrix inequalities in system and control theory},
  volume~15 of {\em SIAM Studies in Applied Mathematics}.
\newblock Society for Industrial and Applied Mathematics (SIAM), Philadelphia,
  PA, 1994.

\bibitem{PDL:Feron}
E.~Feron, P.~Apkarian, and P.~Gahinet.
\newblock Analysis and synthesis of robust control systems via
  parameter-dependent {L}yapunov functions.
\newblock {\em IEEE Trans. Automat. Control}, 41(7):1041--1046, 1996.

\bibitem{PDL:Gahinet}
P.~Gahinet, P.~Apkarian, and M.~Chilali.
\newblock Affine parameter-dependent {L}yapunov functions and real parametric
  uncertainty.
\newblock {\em IEEE Trans. Automat. Control}, 41(3):436--442, 1996.

\bibitem{ParriloGatermann}
K.~Gatermann and P.~A. Parrilo.
\newblock Symmetry groups, semidefinite programs, and sums of squares.
\newblock {\em J. Pure Appl. Algebra}, 192(1-3):95--128, 2004.

\bibitem{PDL:Haddad}
W.~Haddad and D.~S. Bernstein.
\newblock Parameter-dependent {L}yapunov functions and the {P}opov criterion in
  robust analysis and synthesis.
\newblock {\em IEEE Trans. Automat. Control}, 40(3):536--543, 1995.

\bibitem{GarulliHenrion}
D.~Henrion and A.~Garulli, editors.
\newblock {\em Positive polynomials in control}, volume 312 of {\em Lecture
  Notes in Control and Information Sciences}.
\newblock Springer-Verlag, Berlin, 2005.

\bibitem{HolScherer}
C.W.~J Hol and C.W. Scherer.
\newblock A sum-of-squares approach to fixed-order {$H\sb \infty$}-synthesis.
\newblock In {\em Positive polynomials in control}, volume 312 of {\em Lecture
  Notes in Control and Inform. Sci.}, pages 45--71. Springer, 2005.

\bibitem{Khalil}
H.~Khalil.
\newblock {\em Nonlinear Systems}.
\newblock Macmillan, 1992.

\bibitem{KojimaSOS}
M.~Kojima.
\newblock Sums of squares relaxations of polynomial semidefinite programs.
\newblock Research report B-397, Dept. of Mathematical and Computing Sciences,
  Tokyo Institute of Technology, 2003.

\bibitem{KKK}
M.~Krsti\'c, I.~Kanellakopoulos, and P.~Kokotovi\'c.
\newblock {\em Nonlinear and Adaptive Control Design}.
\newblock Wiley, 1995.

\bibitem{Lasserre}
J.~B. Lasserre.
\newblock Global optimization with polynomials and the problem of moments.
\newblock {\em SIAM J. Optim.}, 11(3):796--817, 2001.

\bibitem{OCA}
W.~Lohmiller and J.~J.~E. Slotine.
\newblock On contraction analysis for nonlinear systems.
\newblock {\em Automatica}, 34:683--696, 1998.

\bibitem{Wang-Michel}
A.~N Michel and K.~Wang.
\newblock Robust stability: perturbed systems with perturbed equilibria.
\newblock {\em System and Control Letters}, (21):155--162, 1993.

\bibitem{SOSlyap}
A.~Papachristodoulou and S.~Prajna.
\newblock On the construction of {Lyapunov} functions using the sum of squares
  decomposition.
\newblock {\em Proceedings of IEEE Conference on Decision and Control}, 2002.

\bibitem{SOS_LMI_Pablo}
P.~A. Parrilo.
\newblock On a decomposition of multivariable forms via {LMI} methods.
\newblock {\em Proceedings of the American Control Conference}, 2000.

\bibitem{Parrilo:Phd}
P.~A. Parrilo.
\newblock {\em Structured semidefinite programs and semialgebraic geometry
  methods in robustness and optimization}.
\newblock PhD thesis, California Institute of Technology, May 2000.
\newblock Available at
  \texttt{http://resolver.caltech.edu/CaltechETD:etd-05062004-055516}.

\bibitem{sdprelax}
P.~A. Parrilo.
\newblock Semidefinite programming relaxations for semialgebraic problems.
\newblock {\em Math. Prog.}, 96(2, Ser. B):293--320, 2003.

\bibitem{SOSTOOLS}
S.~Prajna, A.~Papachristodoulou, P.~Seiler, and P.~A. Parrilo.
\newblock {\em {Sum of Squares Optimization Toolbox for MATLAB - User's
  Guide}}.
\newblock \texttt{www.cds.caltech.edu/sostools/sostools.pdf}.

\bibitem{CAND}
S.~Prajna, A.~Papachristodoulou, P.~Seiler, and P.~A. Parrilo.
\newblock {SOSTOOLS}: Control applications and new developments.
\newblock {\em IEEE International Symposium on Computer Aided Control Systems
  Design}, 2004.

\bibitem{Strogatz}
S.~Strogatz.
\newblock {\em Nonlinear Dynamics and Chaos}.
\newblock Addison-Wesley Publishing Company, Reading, MA, 1994.

\bibitem{SeDuMi}
J.~F. Sturm.
\newblock Using {S}e{D}u{M}i 1.02, a {MATLAB} toolbox for optimization over
  symmetric cones.
\newblock {\em Optim. Methods Softw.}, 11/12(1-4):625--653, 1999.

\bibitem{OPC-NO}
W.~Wang and J.~J. Slotine.
\newblock On partial contraction analysis for coupled nonlinear oscillators.
\newblock {\em Biological Cybernetics}, 92(1), 2005.

\bibitem{Zhou}
K.~Zhou, J.~Doyle, and K.~Glover.
\newblock {\em Robust and Optimal Control}.
\newblock Prentice Hall, 1996.

\end{thebibliography}

\end{document}